\documentstyle[12pt]{article}
\begin{document}
\begin{titlepage}
\begin{flushright}
math.QA/0111233
\end{flushright}
\vskip.3in

\begin{center}
{\Large \bf Exchange relations of Level-two vertex operators of
$U_q(\widehat{sl_2})$}
\vskip.3in
{\large Wen-Li Yang$^{a,b}$
\footnote{E-mail:~wlyang@th.physik.uni-bonn.de}
}
\vskip.2in

{\em $~^{a}$ Institute of Modern Physics, Northwest University, Xian
710069, P.R. China\\
$~^{b}$ Physikalisches Institut der Universitat Bonn, Nussallee 12, 53115
Bonn, Germany}
\end{center}

\vskip 2cm
\begin{center}
{\bf Abstract}
\end{center}
We calculate the exchange relations of vertex operators of
$U_q(\widehat{sl_2})$ at  level-two from its bosonic realization. The
corresponding invertibility relation of type I vertex operators is
also studied.


\end{titlepage}


\def\a{\alpha}
\def\b{\beta}
\def\d{\delta}
\def\e{\epsilon}
\def\g{\gamma}
\def\k{\kappa}
\def\l{\lambda}
\def\o{\omega}
\def\t{\theta}
\def\s{\sigma}
\def\D{\Delta}
\def\L{\Lambda}

\def\G{{\cal G}}
\def\Gk{{\cal G}^{(k)}}
\def\R{{\cal R}}
\def\hR{{\hat{\cal R}}}
\def\C{{\bf C}}
\def\P{{\bf P}}

\def\uqa{{U_q(\widehat{sl_2})}}
\def\uqA22{{U_q(A^{(2)}_2)}}


\def\beq{\begin{equation}}
\def\eeq{\end{equation}}
\def\bea{\begin{eqnarray}}
\def\eea{\end{eqnarray}}
\def\ba{\begin{array}}
\def\ea{\end{array}}
\def\no{\nonumber}
\def\lt{\left}
\def\rt{\right}
\newcommand{\bq}{\begin{quote}}
\newcommand{\eq}{\end{quote}}

\newtheorem{Theorem}{Theorem}
\newtheorem{Definition}{Definition}
\newtheorem{Proposition}{Proposition}
\newtheorem{Lemma}[Theorem]{Lemma}
\newtheorem{Corollary}[Theorem]{Corollary}
\newcommand{\proof}[1]{{\bf Proof. }
        #1\begin{flushright}$\Box$\end{flushright}}

\newcommand{\sect}[1]{\setcounter{equation}{0}\section{#1}}
\renewcommand{\theequation}{\thesection.\arabic{equation}}

\sect{Introduction\label{intro}}
Infinite-dimensional highest weight representations and the corresponding vertex
operator \cite{Fre92} of quantum affine (super) algebras are two ingredients of
great importance in algebraic analysis of lattice integrable models
\cite{Jim94,Yan99}.
The exchange relation of the vertex operators and its
invertibility relation play a key role to construct bosonization of the
corresponding lattice models both in the bulk case \cite{Dav93,Koy94,Hou99}
and the
boundary case \cite{Jim95,Fur99,Yan01}.

A powerful approach for studying the highest weight representations and vertex
operators is the bosonization technique \cite{Fre88, Ber89} which allows one to
explicitly construct these objects in terms of the q-deformed free bosonic and
fermionic fields. In this paper, we study the exchange relations
of  vertex operators of $\uqa$ at level-two
\footnote{The exchange relation of the type I vertex operators of the 19-vertex
model at critical regime has been recently  studied by T. Kojima \cite{Koj00}}
 and the invertibility relations from its bosonization which can be realized by
a q-deformed bosonic free field and a fermionic free field \cite{Bou94,Idz94}.

\sect{The quantum affine algebra $\uqa$}
The symmetric Cartan matrix of the affine Lie algebra $\widehat{sl_2}$
is
\bea
(a_{ij})=\left(
\begin{array}{cc}
2&-2\\
-2&2
\end{array}\right)\no
\eea
where $i,j=0,1$. Quantum affine algebra
$\uqa$ is a $q$-analogue of the universal
enveloping algebra of $\widehat{sl_2}$ generated by the Chevalley
generators $\{e_i,f_i,t_i^{\pm 1},d | i=0,1\}$, where $d$ is the
usual derivation operator.
The defining relations are \cite{Cha98}
\bea
& & t_it_j =t_jt_i,\ \  t_id=dt_i, \ \ [d,e_i]=\delta_{i,0}e_i,\ \
[d,f_i]=-\delta_{i,0}f_i,\no\\
& &t_ie_jt_i^{-1}=q^{a_{ij}}e_j,
\ \ t_if_jt_i^{-1}=q^{-a_{ij}}f_j,\no\\
& &[e_i,f_j] =\delta_{ij} \frac{t_i-t_i^{-1}}{q-q^{-1}},\no\\
& & \sum_{r=0}^{1-a_{ij}}(-1)^r
\left[
\begin{array}{c}
1-a_{ij}\\
r
\end{array}
\right]
(e_i)^re_j(e_i)^{1-a_{ij}-r}=0,~~~{\rm if}~i\neq j,\no\\
& & \sum_{r=0}^{1-a_{ij}}(-1)^r
\left[
\begin{array}{c}
1-a_{ij}\\
r
\end{array}
\right]
(f_i)^rf_j(f_i)^{1-a_{ij}-r}=0,~~~{\rm if}~i\neq j,\no
\eea
\noindent where
\bea
& &[n]=\frac{q^n-q^{-n}}{q-q^{-1}},
~~~[n]!=[n][n-1]\cdot\cdot\cdot [1],\no\\
& &\left[
\begin{array}{c}
n\\
r
\end{array}
\right]=\frac{[n]!}{[n-r]!
[r]!}.\label{Note}
\eea

$\uqa$ is a quasi-triangular Hopf algebra
endowed with  Hopf algebra structure:

\begin{eqnarray*}
&&\Delta(t_i)=t_i\otimes t_i,~~~~
   \Delta(e_i)=e_i\otimes 1+t_i\otimes e_i,~~~~
    \Delta(f_i)=f_i\otimes t_i^{-1}+1\otimes f_i,\\
\label{counit}
&&\epsilon(t_i)=1,~~~~\epsilon(e_i)=\epsilon(f_i)=0,\\
&&S(e_i)=-t_i^{-1} e_i,~~~~ S(f_i)=-f_i t_i, ~~~~
   S(t_i^{\pm 1})=t_i^{\mp 1},~~~~ S(d)=-d.
\end{eqnarray*}

$\uqa$ can also be realized by
the Drinfeld generators \cite{Dri88} $\{d,\ \ X^{\pm}_{m}$,
$a_n$, $K^{\pm 1},\gamma^{\pm 1/2} |  m \in
{\bf Z}, n \in {\bf Z}_{\ne 0}\}$.
The relations read
\begin{eqnarray*}
\label{DRB1}
& &\gamma\ \  {\rm is\ \  central },\ \
[K,a_n]=0,\ \ [d,K]=0,\ \ [d,a_n]=na_n,\\
& &[a_m,a_n] =\delta_{m+n,0}
{[2m](\gamma^m-\gamma^{-m})
\over m(q-q^{-1})} ,\\
& &KX^{\pm}_m =q^{\pm 2}X^{\pm}_m K,\ \
[d,X^{\pm}_m]=mX^{\pm}_m ,\\
& &[a_m,X^{\pm}_n]=\pm {[2m] \over m}
\gamma^{\mp |m|/2}X^{\pm}_{n+m},\\
& &[X^{+}_m,X^{-}_n]=\frac{1}{ q-q^{-1}}
(\gamma^{(m-n)/2}\psi^{+}_{m+n} -\gamma^{-(m-n)/2 }
\psi^{-}_{m+n}),\\
\label{DRB2}
& &(z-wq^{\pm 2})X^{\pm}(z)X^{\pm}(w)=
(zq^{\pm 2}-w)X^{\pm}(w)X^{\pm}(z),
\end{eqnarray*}
\noindent where the corresponding Drinfeld currents
$\psi^{\pm}(z)$ and $X^{\pm}(z)$ are defined by
\begin{eqnarray*}
& &\psi^+(z)=\sum_{m=0}^{\infty}\psi^{+}_mz^{-m}=
Kexp\{(q-q^{-1})\sum_{k=1}^{\infty}a_kz^{-k}\},\\
& &\psi^-(z)=\sum_{m=0}^{\infty}\psi^{-}_{-m}z^{m}=
K^{-1}exp\{-(q-q^{-1})\sum_{k=1}^{\infty}a_{-k}z^{k}\},\\
& &X^{\pm}(z)=\sum_{n\in Z}X^{\pm}_mz^{-m-1}.
\end{eqnarray*}

The Chevalley generators are related to the Drinfeld generators
by the formulae:
\bea
\label{CB1}
& &t_1= K,\ \ e_1 = X^{+}_0,\ \ t_0=\gamma K^{-1},
\ \ f_1 = X^{-}_0,\no\\
\label{CB2}
& &e_0= X^{-}_1t_1^{-1},~~~
f_0=t_1X^{+}_{-1}\no
\eea

\subsection{Bosonization of $\uqa$ at level-two}

Let us introduce the $q$-bosonic-oscillators
$\{a_n,~Q,~ P| n\in Z-\{0\}\}$ and  Neuveu Schwartz sector of
$q$-fermionic-oscillators $\{ b_r~|~r\in Z+\frac{1}{2}\}$
which satisfy the commutation relations
\begin{eqnarray*}
& &[a_m,a_n]=\delta_{m+n,0}\frac{[2m]^2}{m},~~m,n\neq 0,
[P,a_m]=[Q,a_m]=0~~,~~[P,Q] = 1,\\
& &\{b_r~,~b_s\}=\frac{[4r]}{2[2r]}\delta_{r+s,0}.
\end{eqnarray*}

\noindent Then  we have \cite{Bou94,Idz94}
\begin{Theorem}
The Drinfeld currents of $\uqa$ at level-two  are realized as
\bea
& & \gamma=q^2,~~ K=q^{2P},\label{CUR1}\\
& &\psi^{+}(z)=q^{2P}exp\{(q-q^{-1})\sum_{n=1}^{\infty}a_nz^{-n}\},\\
& &\psi^{-}(z)=q^{-2P}exp\{-(q-q^{-1})\sum_{n=1}^{\infty}a_{-n}z^{n}\},\\
& &X^{\pm}(z) =\sqrt{2}B(z)E^{\pm}(z), \label{CUR2}
\eea
where
\bea
&&E^{\pm}(z)=
exp\{\pm\sum_{n=1}^{\infty}\frac{a_{-n}}{[2n]}q^{\mp n}z^n\}
exp\{\mp\sum_{n=1}^{\infty}\frac{a_{n}}{[2n]}q^{\mp n}z^{-n}\}
e^{\pm Q}z^{\pm P},\no\\
&&B(z)=\sum_{r\in Z+\frac{1}{2}}b_rz^{-r-\frac{1}{2}}.\no
\eea
\end{Theorem}
\subsection{Level-two vertex operators}
Let $V$ be the 3-dimensional (or spin-$1$) evaluation representation
of $\uqa$, $\{
v_1,v_0,v_{-1}\}$ be the basis vectors of $V$. Then the 3-dimensional
level-0 representation $V_z$ of $\uqa$ is given by \cite{Idz93}
\bea
& &e_1v_m=[1+m]v_{m-1},~~~f_1v_m=[1-m]v_{m+1},
~~~e_0=zf_1,~~~f_0=z^{-1}e_1.
\eea
We  define the dual modules $V_z^*$ of $V_z$ by
$\pi_{V^{*}}(a)=\pi_{V}(S(a))^{t}$, $\forall a\in \uqa$, where
$t$ is the transposition operation.

Throughout, we denote by $V(\l)$ a level-two irreducible
highest weight $\uqa$-module with the highest weight $\l$.
Consider the following intertwines of
$\uqa$-modules:
\begin{eqnarray*}
& &\Phi_{\lambda}^{\mu V}(z) :
 V(\lambda) \longrightarrow V(\mu)\otimes V_{z} ,\ \ \ \
\Phi_{\lambda}^{\mu V^{*}}(z) :
 V(\lambda) \longrightarrow V(\mu)\otimes V_{z}^{*} ,\\
& &\Psi_{\lambda}^{V \mu}(z) :
 V(\lambda) \longrightarrow V_{z}\otimes V(\mu),\ \ \ \
\Psi_{\lambda}^{V^* \mu}(z) :
 V(\lambda) \longrightarrow V_{z}^{*}\otimes V(\mu).
\end{eqnarray*}
They are intertwines in the sense that for any $x\in \uqa$,
\bea
\Theta(z)\cdot x=\Delta(x)\cdot \Theta(z),\ \ \ \Theta(z)=
\Phi(z),\Phi^{*}(z),\Psi(z),\Psi^{*}(z).\no
\eea
$\Phi(z)$ ($\Phi^{*}(z)$) is called type I (dual) vertex operator
and $\Psi(z)$
($\Psi^{*}(z)$) type II (dual) vertex operator.

We expand the vertex operators as
\begin{eqnarray*}
& &\Phi(z)=\sum_{j=1,0,-1}\Phi_j(z)\otimes v_j\  ,\ \ \ \
\Phi^{*}(z)=\sum_{j=1,0,-1}\Phi^{*}_j(z)\otimes v^{*}_j,\\
& &\Psi(z)=\sum_{j=1,0,-1}v_j\otimes\Psi_j(z)\  ,\ \ \ \
\Psi^{*}(z)=\sum_{j=1,0,-1}v^{*}_j\otimes\Psi^{*}_j(z).
\end{eqnarray*}
Define the operators $\phi_j(z),\phi^{*}_j(z),\psi_j(z)$ and
$\psi^{*}_j(z)$ $(j=1,0,-1)$ bosonized by
\bea
&&\phi_1(z)=exp\{\sum_{n=1}^{\infty}\frac{q^{5n}z^n}{[2n]}a_{-n}\}
exp\{-\sum_{n=1}^{\infty}\frac{q^{-3n}z^{-n}}{[2n]}a_n\}~e^Q(-zq^4)^P,
\label{VoxI}\\
&&\phi_0(z)=[\phi_1(z),f_1]_{q^2},
~~\phi_{-1}(z)=\frac{1}{[2]}[\phi_0(z),f_1],
~~\phi^*_i(z)=\phi_{-i}(zq^{-2}),\\
&&\psi_{-1}(z)=exp\{-\sum_{n=1}^{\infty}\frac{q^{n}z^n}{[2n]}a_{-n}\}
exp\{\sum_{n=1}^{\infty}\frac{q^{-3n}z^{-n}}{[2n]}a_n\}~e^{-Q}(-zq^2)^{-P},
\label{VoxII}\\
&&\psi_0(z)=[\psi_{-1}(z),e_1]_{q^2},~~\psi_1(z)=\frac{1}{[2]}[\psi_0(z),e_1],
~~\psi^*(z)=\psi_{-i}(zq^{-2}),
\eea
where $[a,b]_x=ab-xba$. Introduce $\phi(z), \phi^*(z), \psi(z), \psi^*(z)$
by
\bea
& &\phi(z)=\sum_{j=1,0,-1}\phi_j(z)\otimes v_j\  ,\ \ \ \
\phi^{*}(z)=\sum_{j=1,0,-1}\phi^{*}_j(z)\otimes v^{*}_j,\no\\
& &\psi(z)=\sum_{j=1,0,-1}v_j\otimes\psi_j(z)\  ,\ \ \ \
\psi^{*}(z)=\sum_{j=1,0,-1}v^{*}_j\otimes\psi^{*}_j(z).\no
\eea
Using the method of Idzumi \cite{Idz94}, we have
the following result.
\begin{Proposition}
The operators $\phi(z),\phi^*(z), \psi(z), \psi^*(z)$ satisfy the same
commutation relations as  $\Phi^{\mu V}_{\l}(z),\Phi^{\mu
V^*}_{\l}(z),\Psi^{V\mu }_{\l}(z),\Phi^{V^*\mu }_{\l}(z)$ respectively
have, respectively.
\end{Proposition}

To prove the proposition, the relations in appendix B are useful.

\noindent {\bf Remark.} The vertex operators (both type I and type II) can
almost be determined by the method used for the level-one bosonization of
the vertex operators of $U_q(\widehat{sl_N})$ \cite{Koy94}. But the
contribution to the vertex operators of the fermionic part can not be
determined by studying the commutation relations with the
q-bosonic-oscillators. However the commutation relations: $[\Phi_1(z),
X^+(w)]=[\Psi_{-1}(z), X^-(w)]=0$ enables us to determine them  uniquely
up to some scalar factor.
\sect{Level-two highest weight $\uqa$-modules and the corresponding
intertwines}
Set $ {\cal F}_1=\oplus_{n\in
Z}C[a_{-1},a_{-2},.....;~b_{-\frac{1}{2}},b_{-\frac{3}{2}},....]
e^{nQ}|0>$, where the Fock vacuum vector $|0>$ is defined by
\begin{eqnarray*}
&&a_n|0>=0,~~~{\rm for}~ n>0,~~~~P|0>=0,\\
&&b_{l+\frac{1}{2}}|0>=0~~{\rm for}~l\geq 0.
\end{eqnarray*}

It can be shown that the bosonized action of $\uqa$ on ${\cal F}_1$ is
closed. Hence the Fock space constitutes a $\uqa$-module at level-two.
However, it is not irreducible. In order to obtain the irreducible
subspace in ${\cal F}_1$, we should introduce the GSO-like projectors
\cite{Bou94}:
\bea
P_{\pm}=\frac{1\pm exp\{-2\pi i~d\}}{2},
\eea
where
\bea
d=-\sum_{n=1}^{\infty}\frac{n^2}{[2n]^2}a_{-n}a_n-
\sum_{l=0}^{\infty}
\frac{(2l+1)[2l+1]}{[4l+2]}b_{-l-\frac{1}{2}}b_{l+\frac{1}{2}}
-\frac{P^2}{2}.\no
\eea
Note that $[d,a_n]=na_n,~[d,b_r]=rb_r,~[d,Q]=P$, we have the following
relations
\bea
&&x^{-d}\phi_1(z)x^{d}=x^{\frac{1}{2}}\phi_1(xz),~~
x^{-d}\psi_{-1}(z)x^{d}=x^{\frac{1}{2}}\psi_{-1}(xz),\\
&&x^{-d}X^{\pm}(z)x^{d}=xX^{\pm}(xz).
\eea
Then we have $P_{\pm}X=XP_{\pm}$ for any $X\in \uqa$.
Define
${\cal F}^{(0)}=P_{+}{\cal F}_1$ and ${\cal F}^{(1)}=P_{-}{\cal F}_1$.
Then we have \cite{Bou94}
\begin{Theorem}
${\cal F}^{(0)}$ and ${\cal F}^{(1)}$ are the irreducible highest weight
$\uqa$-modules with the highest weights $2\L_0$ and $2\L_1$ respectively,
namely, ${\cal F}^{(0)}=V(2\L_0)$ and ${\cal F}^{(1)}=V(2\L_1)$.
\end{Theorem}

One can further find that $P_{\pm}\Theta (z)=\Theta (z) P_{\mp}$
for $\Theta (z)=\phi_i(z),\phi^*_i(z),\psi_i(z),\psi^*_i(z)$.
Hence, we have following equality:
\bea
&&\Phi^{2\L_1~V}_{2\L_0}(z)=P_-\phi(z) P_+,~~
\Phi^{2\L_0~V}_{2\L_1}(z)=P_+\phi(z) P_-,\label{VoxI1}\\
&&\Psi^{V~2\L_1}_{~2\L_0}(z)=P_-\psi(z) P_+,~~
\Psi^{V~2\L_0}_{~2\L_1}(z)=P_+\psi(z) P_-,\label{VoxII1}\\
&&P_{\pm}\Theta(z)P_{\pm}=0,~~{\rm for}~
\Theta (z)=\phi_i(z),\phi^*_i(z),\psi_i(z),\psi^*_i(z).\label{Pro}
\eea
\sect{Exchange relations of vertex operators}
In this section, we derive the exchange relations of the type I and type
II
operators of $\uqa$ at level-two from their bosonization.
\subsection{The R-matrix}
Let $R(z) \in End(V\otimes V)$ be the R-matrix of $\uqa$ defined by
\bea
R(z)(v_i\otimes v_j)=\sum_{k,l}R^{ij}_{kl}(z) v_k\otimes v_l.\no
\eea
It can be given explicitly by
\bea
R(z)=r(z)\left(
\begin{array}{ccccccccc}
1 & 0 & 0 & 0 & 0 & 0 & 0 & 0 & 0\\
0 & b & 0 & e & 0 & 0 & 0 & 0 & 0\\
0 & 0 & d & 0 & g & 0 & f & 0 & 0\\
0 & \bar{e} & 0 & b & 0 & 0 & 0 & 0 & 0\\
0 & 0 & \bar{g} & 0 & a & 0 & h & 0 & 0\\
0 & 0 & 0 & 0 & 0 & b & 0 & e & 0\\
0 & 0 & \bar{f} & 0 & \bar{h} & 0 & d & 0 & 0\\
0 & 0 & 0 & 0 & 0 & \bar{e} & 0 & b & 0\\
0 & 0 & 0 & 0 & 0 & 0 & 0 & 0 & 1
\end{array}
\right). \label{r12}
\eea
Here the normalized partition is
\bea
&&r(z)=\frac{1-zq^2}{z-q^2},
\eea
which is just that of Ref.\cite{Idz93} with the level $k=2$. The other
nonzero elements are given by
\bea
&&a=\frac{qz^2-(2q+2q^{-1}-q^3-q^{-3})z+q^{-1}}{(zq^2-q^{-2})(zq-q^{-1})},\\
&&b=\frac{z-1}{zq^2-q^{-2}},~~~
d=\frac{(z-1)(zq^{-1}-q)}{(zq^2-q^{-2})(zq-q^{-1})},\\
&&e=\frac{q^2-q^{-2}}{zq^2-q^{-2}},
~~~\bar{e}=\frac{(q^2-q^{-2})z}{zq^2-q^{-2}},\\
&&f=\frac{(q^2-q^{-2})(qz^2-(q+q^{-1})z+q)}
{(zq^2-q^{-2})(zq-q^{-1})},\\
&&\bar{f}=\frac{(q^2-q^{-2})(q^{-3}z-q^{-3}+q-q^{-1})z}
{(zq^2-q^{-2})(zq-q^{-1})},\\
&&g=\frac{(q^2-q^{-2})(q^2+1)(z-1)}
{(zq^2-q^{-2})(zq-q^{-1})},\\
&&\bar{h}=\frac{(q^2-q^{-2})(q^2+1)(z-1)z}
{(zq^2-q^{-2})(zq-q^{-1})},\\
&&\bar{g}=\frac{q^{-2}}{[2]^2}~\bar{h},~~~h=\frac{q^{-2}}{[2]^2}~g.
\eea
The R-matrix satisfies the Yang-Baxter equation on $V\otimes V\otimes V$
\bea
R_{12}(z)R_{13}(zw)R_{23}(w)=R_{23}(w)R_{13}(zw)R_{12}(z),\no
\eea
and the initial condition $R(1)=P$ with $P$ being the permutation
operator.
\subsection{The exchange relations}
Define
\bea
\oint dz f(z)=Res(f)=f_{-1},~~ {\rm for~ formal~series~function~}
f(z)=\sum_{n\in Z}f_nz^n.
\eea
Then the chevalley generators of $\uqa$ can be expressed by the integrals
\bea
e_1=\oint dz X^{+}(z),~~f_1=\oint dz X^{-}(z).
\eea
From the normal order relations in appendix A, one can also obtain the
integral expression of the vertex operators defined in
(\ref{VoxI}-\ref{VoxII})
\bea
&&\phi_0(z)=\sqrt{2}\oint dw \frac{(q^2-q^{-2})}{(zq^2)(1-\frac{w}{zq^2})
(1-\frac{zq^6}{w})}:\phi_1(z)E^-(w):B(w),\\
&&\phi_{-1}(z)=\frac{2}{[2]}\oint dw \oint d\eta~
\frac{(q^2-q^{-2})}{(zq^2)(1-\frac{w}{zq^2})
(1-\frac{zq^6}{w})}:\phi_1(z)E^-(w)E^-(\eta):\no\\
&&~~~~~~~~~~~~~\times \lt\{\frac{w-\eta
q^2}{(-zq^4)(1-\frac{\eta}{zq^2})}B(w)B(\eta)~-~
\frac{\eta-wq^2}{\eta(1-\frac{zq^6}{\eta})}B(\eta)B(w)\rt\},\\
&&\psi_0(z)=\sqrt{2}\oint dw \frac{(q^{-2}-q^{2})}
{(1-\frac{w}{zq^4})w
(1-\frac{z}{w})}:\psi_{-1}(z)E^+(w):B(w),\\
&&\psi_{1}(z)=\frac{2}{[2]}\oint dw \oint d\eta~
\frac{(q^{-2}-q^{2})}{(1-\frac{w}{zq^4})w
(1-\frac{z}{w})}:\phi_1(z)E^+(w)E^+(\eta):\no\\
&&~~~~~~~~~~~~~\times \lt\{\frac{-w(1-\frac{\eta}{wq^2})
}{(zq^2)(1-\frac{\eta}{zq^4})}B(w)B(\eta)~-~
\frac{\eta(1-\frac{w}{\eta q^2})
}{\eta(1-\frac{z}{\eta})}B(\eta)B(w)\rt\}.
\eea
By the ``weak equality" technique proposed in \cite{Asa96}, using the
normal order relations given in appendix A and the relations in appendix
B, after tedious calculations , we can show that the bosonic vertex
operators defined in (\ref{VoxI}-\ref{VoxII}) satisfy the
Faddeev-Zamolodchikov (ZF) algebra
\bea
& &\phi_j(z_2)\phi_i(z_1)=\sum_{kl}R(\frac{z_1}{z_2})^{kl}_{ij}
\phi_k(z_1)\phi_l(z_2),\no\\
& &\psi_i(z_1)\psi_j(z_2)=\sum_{kl}R(\frac{z_1}{z_2})^{kl}_{ij}
\psi_l(z_2)\psi_k(z_1),\no\\
& &\psi_i(z_1)\phi_j(z_2)=\tau(\frac{z_1}{z_2})\phi_j(z_2)\psi_i(z_1),\no
\eea
where $\tau(z)=-1$. Furthermore, we can show that the type I bosonic
vertex operators have the following invertibility relation
\bea
&& \phi_i(z)\phi^*_i(z)=f_i~id,\no
\eea
where the scalar factors $\{f_i\}$ are $f_1=\frac{1}{z(q^2-1)q^2}$,
$f_0=-[2]\frac{1}{z(q^2-1)q^2}$,  $f_{-1}=\frac{q^{-2}}{z(q^2-1)q^2}$.
In the derivation of the above relations we have used the identity
$:\phi_1(z)\phi_1(zq^{-2})E^-(zq^4)E^-(zq^2):=id$.

From the theorem 2 and the properties of bosonized vertex operators
(\ref{Pro}), we have
\begin{Proposition}
The vertex between the irreducible level-two highest weight $\uqa$-module
$V(2\L_0)$ and $V(2\L_1)$ satisfy the ZF algebra and invertibility
relation
\bea
& &\Phi_j(z_2)\Phi_i(z_1)=\sum_{kl}R(\frac{z_1}{z_2})^{kl}_{ij}
\Phi_k(z_1)\Phi_l(z_2),\\
& &\Psi_i(z_1)\Psi_j(z_2)=\sum_{kl}R(\frac{z_1}{z_2})^{kl}_{ij}
\Psi_l(z_2)\Psi_k(z_1),\\
& &\Psi_i(z_1)\Phi_j(z_2)=\tau(\frac{z_1}{z_2})\Phi_j(z_2)\Psi_i(z_1),\\
&& \Phi_i(z)\Phi^*_i(z)=f_i~id.
\eea
\end{Proposition}

This agrees with the results of \cite{Idz93} which were obtained from the
solution of
Q-KZ equation.

\noindent {\bf Remark.} We can similarily derive that the level-two vertex
operators among $V(\L_0+\L_1)$s with the 3-dimensional evaluation representation
in Ref.\cite{Idz94} satisfies the  ZF algebraic relation with the same
R-matrix.

\vskip.3in
\noindent {\bf Acknowledgments.}
The author   would like to thank
Prof. G. von Gehlen for his encouragements and useful comments.
This work has
been supported by
the Alexander von Humboldt Foundation.

\section*{Appendix A}
In this appendix, we give the normal order relations of the fundamental
bosonic fields including the ferimonic field:
\bea
&&\phi_1(z)\phi_1(w)=(-zq^4)(1-\frac{wq^2}{z}):\phi_1(z)\phi_1(w):,\no\\
&&\phi_1(z)E^+(w)=(-zq^4)(1-\frac{w}{zq^4}):\phi_1(z)E^+(w):
=E^+(w)\phi_1(z),\no\\
&&\phi_1(z)E^-(w)=\frac{1}{(-zq^4)(1-\frac{w}{zq^2})}:\phi_1(z)E^-(w):,\no\\
&&E^-(w)\phi_1(z)=\frac{1}
{w(1-\frac{zq^6}{w})}:\phi_1(z)E^-(w):,\no\\
&&\psi_{-1}(z)\psi_{-1}(w)=(-zq^2)(1-\frac{w}{zq^2}):\psi_{-1}(z)
\psi_{-1}(w):,\no\\
&&\psi_{-1}(z)E^-(w)=(w-zq^2):\psi_{-1}(z)E^-(w):
=E^-(w)\psi_{-1}(z),\no\\
&&\psi_{-1}(z)E^+(w)=-\frac{1}{(zq^2)(1-\frac{w}{zq^4})}:\psi_{-1}(z)E^+(w):,
\no\\
&&E^+(w)\psi_{-1}(z)=\frac{1}
{w(1-\frac{z}{w})}:\psi_{-1}(z)E^+(w):,\no\\
&&E^+(z)E^+(w)=(z-wq^{-2}):E^+(z)E^+(w):,\no\\
&&E^-(z)E^-(w)=(z-wq^{2}):E^-(z)E^-(w):,\no\\
&&E^+(z)E^-(w)=\frac{1}{z(1-\frac{w}{z})}:E^+(z)E^-(w):,\no\\
&&E^-(w)E^+(z)=\frac{1}{w(1-\frac{z}{w})}:E^-(w)E^+(z):,\no\\
&&B(z)B(w)=\frac{[2]}{2}\frac{1-\frac{w}{z}}{z(1-\frac{wq^2}{z})
(1-\frac{w}{zq^{2}})}+:B(z)B(w):.\no
\eea

\section*{Appendix B}
By means of the bosonic realization of $\uqa$, the integral expressions of
the vertex operators and the ``weak equality" technique given in
Ref.\cite{Asa96}, one can
check the following relations.

\begin{itemize}
\item For the type I vertex operators
\bea
&& [\phi_1(z), e_1]=0,~~[\phi_0(z),e_1]=[2]t_1\phi_1(z),\no\\
&&[\phi_{-1}(z),e_1]=t_1\phi_0(z),~~[\phi_{-1}(z),f_1]_{q^{-2}}=0,\no\\
&&t_1\phi_1(z)=q^2\phi_1(z)t_1,~~t_1\phi_0(z)=\phi_0(z)t_1,
~~t_1\phi_{-1}(z)=q^{-2}\phi_{-1}(z)t_1.\no
\eea

\item For the type II vertex operators
\bea
&& [\psi_{-1}(z), f_1]=0,~~[\psi_0(z),f_1]=[2]t^{-1}_1\psi_{-1}(z),\no\\
&&[\psi_{1}(z),f_1]=t^{-1}_1\psi_0(z),
~~[\psi_{1}(z),e_1]_{q^{-2}}=0,
\no\\
&&t_1\psi_1(z)=q^2\psi_1(z)t_1,~~t_1\psi_0(z)=\psi_0(z)t_1,
~~t_1\psi_{-1}(z)=q^{-2}\psi_{-1}(z)t_1.\no
\eea

\end{itemize}

\bibliographystyle{unsrt}

\end{document}